 \input amstex
\documentstyle{amsppt}

\magnification=\magstep1
\parindent=1em 
     \baselineskip=18pt      

\hsize=16 true cm
\vsize=24 true cm

 %
    %
           %
\def\ove#1{\overline{#1}}    %

     \def\({\left(}       \def\al{\alpha}           
     \def\){\right)}      \def\e{\varepsilon}    
     \def\[{\left[}       
     \def\]{\right]}      \def\ffi{\varphi}

     \def\<{\langle}                 
     \def\>{\rangle}

        \def\co{\operatorname{co}\,}




             \def\bigp{\bigpagebreak}

\CenteredTagsOnSplits                        \NoBlackBoxes



\NoRunningHeads
\pageno=1  %
\footline={\hss\tenrm\folio\hss} %
    %

        \topmatter
        \title {On dentability in locally convex vector spaces}
        \endtitle
          \author {Oleg Reinov and Asfand Fahad}  \endauthor
\address\newline
          \phantom{AA}              \tenrm {Oleg Reinov:}\newline
St. Petersburg State University, Dept. Math. and Mech.,\newline
198904, Saint Petersburg,
Russia and\newline
Abdus Salam School of Mathematical Sciences, \newline
68-B, New Muslim Town, Lahore 54600, PAKISTAN
\endaddress

\email
                          \tenrm
orein51\@mail.ru
\endemail

\address\newline
     \phantom{AA}                         \tenrm {Asfand Fahad:}  \newline
Abdus Salam School of Mathematical Sciences, \newline
68-B, New Muslim Town, Lahore 54600, PAKISTAN
\endaddress

\email
                          \tenrm
asfandfahad1\@yahoo.com
\endemail


\abstract\nofrills
\endabstract
    \endtopmatter

\document
\footnote""{
The research is supported by the Higher  Education Commission of Pakistan
and by grant 12-01-00216 of RFBR.}

\footnote""{${ }^\maltese$ AMS Subject Classification 2010: 46B22  Radon-Nikodym, Krein-Milman and related properties;
46A55 Convex sets in topological linear spaces; Choquet theory.}
\footnote""{${ }$ Key words: dentability, dentable sets, locally convex spaces. }

 \bigp

 \head \S 0. Preliminaries\endhead

For a locally convex vector space (l.c.v.s.) $E$ and an (absolutely convex) neighbor\-hood $V$ of zero,
a bounded subset $A$ of $E$ is said to be {\it $V$-dentable}\ (respectively, {\it $V$-f-dentable}\,)
 if for any $\epsilon>0$ there exists an $x\in A$ so that
 $$
 x\notin \overline{\co}\, (A\setminus (x+\epsilon V))
 $$
 (respectively, so that
  $$
 x\notin {\co}\, (A\setminus (x+\epsilon V))\, ).
 $$
 Here, "$\overline{\co}$" denotes  the closure in $E$ of the convex hull of a set.

 We present a theorem which says that for a wide class of bounded subsets $B$
 of locally convex vector spaces the following is true:
 \smallskip

 $(V)$\ every subset of $B$ is $V$-dentable if and only if
 every subset of $B$ is $V$-f-dentable.
 \smallskip

The proof is purely geometrical and independent of any related facts.

As a consequence (in the particular case where $B$ is complete convex bounded metrizable
subset of a l.c.v.s.), we obtain a positive solution to a 1978-hypothesis of an
American mathematician Elias Saab (see p. 290 in
"On the Radon-Nikodym property in a
class of locally convex spaces",
Pacific J. Math. 75, No. 1, 1978,  281-291).

\smallskip


 \head \S 1. A proposition\endhead

\proclaim {\bf  Proposition.}\it
Let $B$ be a bounded sequentially complete convex metrizable subset of
a locally convex vector space $E,$ $V$ is a neighborhood of zero in $E.$
The following are equivalent:

$1).$\, $B$ is subset $V$-dentable\footnote{"subset $\Pi$" means that every subset of a set is $\Pi.$};

$2).$\, $B$ is subset $V$-f-dentable.
\endproclaim\rm

   \demo{\it Proof}\
 Clearly, $1) \implies 2).$
 Let $d$ be a metric in $B$ which gives the topology on $B,$ induced from $E.$
 For $x\in B,$ denote by $\ove{D}_\e(x)$ the $e$-ball in $B$ with a center at $x.$
 We may and do assume that $V$ is absolutely convex.

 Given 2), suppose that there is a subset $K\subset B,$ which is not $V$-dentable, i.e.
 there exists an $\e>0$ such that $x\in \ove{\co}(K\setminus (2\e V+x))$ for all $x\in K.$
 Put (for a later use) $k(i,0)=0$ if $i=1,2,\dots.$ Fix $x_{01}\in K$ and take $x_{1j}\in K$
$ (j=1,2, \dots, k(0,1)$) such that
$$
 x_{1j}\notin x_{01} + 2\e V,\ \text{ i.e. }\ x_{1j}-x_{01}\notin 2\e V;
$$
 $$
  d(x_{01}, \sum_{j=1}^{k(0,1)} \al_{1j}x_{1j}) <\e_0/8,
 $$
 and
 $$
  \ove{D}_{\e_0}(x_{01})\subset \e/8\, V+x_{01},
 $$
 where $\al_{1j}\ge0,$ $\sum_j \al_{1j}=1$ and $\e_0$ is such that $\e_0<\e/8.$

 Let $\e_1>0$ be such a number that
 $0<\e_1<\e_0/2,$  $x_{1j}-x_{01}\notin (\e+\e_1) V,$ $\ove{D}_{\e_1(x_{1j})}\subset \e/8\, V+x_{1j}$
 and $\ove{D}_{\e_1}(x_{1j})\cap [(\e+\e_1) V+x_{01}]=\emptyset$
  for all $j=1,2,\dots, k(0,1)=: n(0).$

 For every $x_{1m}$ take $(x_{2j})\subset K$
 $(j= k(1,m-1)+1,\dots, k(1,m))$\,
 so that
 $$
 x_{2j}- x_{1m}\notin 2\e V,
 $$
 $$
  d(x_{1m}, \sum_{j=k(1,m-1)+1}^{k(1,m)} \al_{2j} x_{2j}) < \e_1/4
 $$
 and
 $$
d( \sum_{j_1=1}^{k(0,1)} \al_{1j_1} \sum_{j_2=k(1,m-1)+1}^{k(1,m)} \al_{2j_2} x_{2j_2},  \sum_{j_1=1}^{k(0,1)} \al_{1j_1} x_{1j_1})<\e_1/4,
 $$
 where $\al_{2j}\ge0,$ $\sum_{j_2=k(1,m-1)+1}^{k(1,m)} \al_{2j_2}=1$
 (we use the continuity of the metric $d;$ see the next step).

 Let $n(1)=k(1,n(0))$ and $\e_2$ be such that
 $0<\e_2<\e_1/2,$  $x_{1m}-x_{2j}\notin (\e+\e_2) V,$ $\ove{D}_{\e_2(x_{2j})}\subset \e/8\, V+x_{2j}$
 and $\ove{D}_{\e_2}(x_{2j})\cap [(\e+\e_2) V+x_{1m}]=\emptyset$
  for every pair of indices $m=1,2,\dots, n(0)$ and
  $j=k(1,m-1)+1,\dots, k(1,m).$

 Now, find  a $\delta_2>0, \delta_2<\e_2/4^2,$ with a property that it follows from $(\tilde{x}_{2j})\subset B$ and
 $\max_j \, d(x_{2j}, \tilde{x}_{2j})<\delta_2$ that for all
 $m=1,2,\dots, k(0,1)$
 $$
  d( \sum_{j_2=k(1,m-1)+1}^{k(1,m)} \al_{2j_2} x_{2j_2},  \sum_{j_2=k(1,m-1)+1}^{k(1,m)} \al_{2j_2} \tilde{x}_{2j_2})<\e_2/4^2
 $$
 and
 $$
 d( \sum_{j_1=1}^{k(0,1)} \al_{1j_1} \sum_{j_2=k(1,m-1)+1}^{k(1,m)} \al_{2j_2} x_{2j_2},
 \sum_{j_1=1}^{k(0,1)} \al_{1j_1} \sum_{j_2=k(1,m-1)+1}^{k(1,m)} \al_{2j_2} \tilde{x}_{2j_2})<\e_2/4^2
 $$

 For every $x_{2m}$ take $(x_{3j})\subset K$
 $(j=k(2,m-1)+1, \dots, k(2,m))$ so that
 $$
 x_{2m}- x_{3j}\notin 2\e V,
 $$
 $$
  d(x_{2m}, \sum_{j=k(2,m-1)+1}^{k(2,m)} \al_{3j} x_{3j}) <\delta_2 \, (<\e_2/4^2),
 $$
 where $\al_{3j}\ge0,$ $\sum_{j=k(2,m-1)+1}^{k(2,m)} \al_{3j}=1.$
 Putting $\tilde{x}_{2j}:=  \sum_{j=k(2,m-1)+1}^{k(2,m)} \al_{3j} x_{3j},$
 we get the corresponding three inequalities combining the last ones.

 And one more step (of induction).

 Let $n(2)=k(1,n(1))$ and let $\e_3$ be such that
  $0<\e_3<\e_2/2,$  $x_{2m}-x_{3j}\notin (\e+\e_3) V,$ $\ove{D}_{\e_3(x_{3j})}\subset \e/8\, V+x_{3j}$
 and $\ove{D}_{\e_3}(x_{3j})\cap [(\e+\e_3) V+x_{2m}]=\emptyset$
  for every pair of indices $m=1,2,\dots, n(1)$ and
  $j=k(2,m-1)+1,\dots, k(2,m).$

 Find
 a $\delta_3>0, \delta_3<\e_3/4^3$ such that it follows from $(\tilde{x}_{3j})\subset B$ and
 $\max_j \, d(x_{3j}, \tilde{x}_{3j})<\delta_3$ that for all  $m_2=1,2,\dots, n(1)$ and
 $m_1=1,2,\dots, n(0):$
 $$
  d( \sum_{j_3=k(2,m_2-1)+1}^{k(2,m_2)} \al_{3j_3} x_{3j_3},  \sum_{j_3=k(2,m_2-1)+1}^{k(2,m_2)} \al_{3j_3} \tilde{x}_{3j_3})<\e_3/4^3,
 $$
 $$  \multline
 d( \sum_{j_2=k(1,m_1-1)+1}^{k(1,m_1)} \al_{2j_2} \sum_{j_3=k(2,m_2-1)+1}^{k(2,m_2)} \al_{3j_3} x_{3j_3},  \\
\sum_{j_2=k(1,m_1-1)+1}^{k(1,m_1)} \al_{2j_2} \sum_{j_3=k(2,m_2-1)+1}^{k(2,m_2)} \al_{3j_3}\tilde{x}_{3j_3})<\e_3/4^3
\endmultline
 $$
 and
 $$ \multline
 d( \sum_{j_1=1}^{k(0,1)} \al_{1j_1} \sum_{j_2=k(1,m_1-1)+1}^{k(1,m_1)} \al_{2j_2} \sum_{j_3=k(2,m_2-1)+1}^{k(2,m_2)} \al_{3j_3} x_{3j_3},  \\
\sum_{j_1=1}^{k(0,1)} \al_{1j_1}\sum_{j_2=k(1,m_1-1)+1}^{k(1,m_1)} \al_{2j_2} \sum_{j_3=k(2,m_2-1)+1}^{k(2,m_2)} \al_{3j_3}\tilde{x}_{3j_3})<\e_3/4^3.
 \endmultline
 $$

For every $x_{3m_2}$
 \, $(m_2=1,\dots, n(1))$ take $(x_{4j})\subset K$\,
 $(j= k(3,m_2-1)+1, \dots, k(3,m_2))$ so that
  $$
 x_{3m_2}- x_{4j}\notin 2\e V,
 $$
 $$
  d(x_{3m_2}, \sum_{j=k(3,m_2-1)+1}^{k(3,m_2)} \al_{4j} x_{4j}) <\delta_3 \, (<\e_3/4^3),
 $$
 where $\al_{4j}\ge0,$ $\sum_{j=k(3,m_2-1)+1}^{k(3,m_2)} \al_{4j}=1.$
 Putting $\tilde{x}_{3j_3}:=  \sum_{j=k(3,m_2-1)+1}^{k(3,m_2)} \al_{4j} x_{4j},$
 we get four corresponding  inequalities combining the last ones.

 Then, let $n(3)=k(2,n(2))$ and let $\e_4$ be such that
  $0<\e_4<\e_3/2,$  $x_{3m}-x_{4j}\notin (\e+\e_4) V,$ $\ove{D}_{\e_4(x_{4j})}\subset \e/8\, V+x_{4j}$
 and $\ove{D}_{\e_4}(x_{4j})\cap [(\e+\e_4) V+x_{3m}]=\emptyset$
  for every corresponding pair $x_{3m}, x_{4j}$ etc.

 We think that induction is clear.

 By this induction, we get a subset $(x_{im})_{i,m}\subset K$
 with properties like ones just indicated above for
 a part of the set. Having this sequences $(x_{im})_{i=0,m=1}^{\infty,n(i-1)}$
 (with $n(-1)=1$) and $(\e_i)_{i=0}^\infty$ in mind, we, for any $q=1,2,\dots;$
 $i=0,1,2,\dots $ and $m=1,2,\dots, n(i-1),$ put
 $$ \multline
  \ffi(i,m;q):=\sum_{j_{m+1}=k(i,m-1)+1}^{k(i,m)} \al_{i+1j_{m+1}}
  \sum_{j_{m+2}=k(i+1, j_{m+1}-1)+1}^{k(i+1,j_{m+1})} \al_{i+2j_{m+2}}  \\
  \sum_{j_{m+3}=k(i+2, j_{m+2}-1)+1}^{k(i+2,j_{m+2})} \al_{i+3j_{m+3}}\dots
  \sum_{j_{m+q}=k(i+q-1, j_{m+q-1}-1)+1}^{k(i+q-1,j_{m+q-1})} \al_{i+qj_{m+q}} x_{i+qj_{m+q}},
  \endmultline
 $$

Then, by triangle inequality, $d(\ffi(i,m;q), x_{im})\le \e_i$ and
$d(\ffi(i,m;q_1), \ffi(i,m;q_2))\underset{q_1,q_2}\to{\to}0.$
Therefore, for any $i$ and $m,$ there exists $z_{im}\in \ove{\co} K$ for which
$\ffi(i,m;q)\underset{q}\to{\to} z_{im}$ in $B$ and

$$(i)\, \ d(z_{im}, x_{im})\le\e_i,
$$
$$
 (ii)\, \ z_{im}=\sum_{j=k(i,m-1)+1}^{k(i,m)} \al_{i+1j} z_{i+1j}.
$$

To obtain the second equality, it is enough to compare the definition of
$\ffi(i,m;q)$ with corresponding formulas for $\ffi(i,m;m+1)$ and $\ffi(i+1,j_{m+1};q)$
to note that
$$
 \ffi(i,m;q+1)=\sum_{j_{m+1}=k(i,m-1)+1}^{k(i,m)} \al_{i+1j_{m+1}} [\ffi(i+1, j_{m+1};q)].
$$
Taking a limit, as $q\to +\infty,$ we get (ii).

Fix $i=1,2,\dots;$ $m=1,2,\dots, n(i-1)$ and
$j=k(i,m-1)+1,\dots, k(i,m).$
Since $x_{im}-x_{i+1j}\notin (\e+\e_{i+1}) V,$
we get from (i) and the construction that, e.g., $z_{im}-z_{i+1j}\notin (\e/10) V.$
Together with (ii), this gives us a contradiction.

  \enddemo

 \head \S 2. Main Theorem\endhead

 The proof of the following simple assertion can be given by a reader:

 \proclaim {\bf  Lemma.}\it
  Let $E$ be a locally convex  vector space and let $ B \subset E $ be closed bounded convex sequentially complete and
having the property that
for every $ M \subset B $ and for $ x\in \overline{M}$ there exists a sequence ${x_{n}}\in M$ such that $\lim \limits_{n} {x_{n}}=x.$
If $B$ is not $V$-dentable then there exists a countable set $A \subseteq B $ which is not $V$-dentable
\endproclaim\rm

Proposition and Lemma yields the following theorem.

\proclaim {\bf  Theorem.}\it
 Let $E$ be a locally convex  vector space, $V$ is a neighborhood of zero in $E$  and let $ B \subset E $ have the  following properties:
  $(i)$\, it is closed bounded convex and sequentially complete;
$(ii)$\, for every $ M \subset B $ and for $ x\in \overline{M}$ there exists a sequence ${x_{n}}\in M$ such that $\lim_{n} {x_{n}}=x;$
$(iii)$\, each separable subset of $B$ is metrizable.

Then the following are equivalent:

  $(i)$\, $B$ is subset $V$-dentable;

 $(ii)$\, $B$ is subset $V$-f-dentable.
 \endproclaim\rm
\enddocument

\bigp


\enddocument